\documentclass[preprint,12pt,authoryear]{elsarticle}

\usepackage[frenchb, english]{babel} 
\usepackage{amssymb}
\usepackage{amsfonts,amsmath,amssymb,amsthm,amsopn,amscd} 
\usepackage{dsfont} 
\usepackage[all]{xy} 
\usepackage{stmaryrd}
\usepackage{float}
\usepackage{here}
\usepackage{cases}

\usepackage{lscape}

\usepackage{color}

\usepackage[T1]{fontenc}
\usepackage{caption}
\usepackage{array,multirow,makecell}
\setcellgapes{1pt}
\makegapedcells
\newcolumntype{R}[1]{>{\raggedleft\arraybackslash }b{#1}}
\newcolumntype{L}[1]{>{\raggedright\arraybackslash }b{#1}}
\newcolumntype{C}[1]{>{\centering\arraybackslash }b{#1}}

\usepackage[utf8x]{inputenc}

\definecolor{gris25}{gray}{0.75}
\definecolor{gris10}{gray}{0.90}

\journal{Elsevier}

\begin{document}

\begin{frontmatter}



\title{EM estimation of a Structural Equation Model}

\author[rvt,els]{X.~Bry}
\ead{xavier.bry@univ-montp2.fr}
\author[rvt,focal]{C.~Lavergne}
\ead{christian.lavergne@univ-montp3.fr}
\author[rvt,els]{M.~Tami\corref{cor1}}
\ead{myriam.tami@univ-montp2.fr}


\cortext[cor1]{Corresponding author}
\address[rvt]{Institut Montpelliérain Alexander Grothendieck (IMAG), CNRS, France}

\address[els]{University of Montpellier, France}
\address[focal]{University Paul-Valéry Montpellier 3, France}

\begin{abstract}
In this work, we propose a new estimation method of a Structural Equation Model. Our method is based on the EM likelihood-maximization algorithm. We show that this method provides estimators, not only of the coefficients of the model, but also of its latent factors. Through a simulation study, we investigate how fast and accurate the method is, and then apply it to real environmental data.
\end{abstract}


\end{frontmatter}


\section{Introduction}
Structural Equation Models (SEM) are widely used is as various research domains as psychology, social and behavioral sciences, ecology, chemometrics, etc. A SEM formalizes the interdependence of many Observed numeric Variables (OV) through fewer unobserved ones, referred to as Latent Variables (LV). Every LV is assumed to be underlying a specific set of OVs, which depend on it as well as on extra observed covariates. A SEM is structured through two types of equations, termed measurement equations and structural equations. A measurement equation relates a LV to the corresponding OV's. A structural equation states a hypothesized relationship between LV's. \ref{schema_p} graphs an example of SEM.\\
Literature widely presents two competing families of methods that deal with SEM's: factor-methods, and component-methods. Among the former family are the classical Factor Analysis, and Jöreskog's SEM estimation technique \cite{joreskog_general_1970} implemented in the LISREL software. These methods use factors as LV's, i.e. variables of which we merely assume to know the distribution (typically standardized normal). They base their estimation on the structure of the covariance matrix of the data according to the model, within a likelihood maximization approach. They estimate all coefficients in the model (linear combination coefficients and variances), but not the values of the factors, which therefore remain unknown. The component-model family of methods assumes that every LV is a component, i.e. a linear combination, of its OV's. Note that such a constraint is stronger than the distribution-assumption made on factors. This family includes the classical Principal Component Analysis, Canonical Correlation Analysis and Partial Least Squares (PLS) methods \cite{joreskog_recent_1982}, \cite{wangen_multiblock_1989}, \cite{lohmoller_latent_2013}, \cite{w._w._chin_structural_1999}, \cite{vinzi_handbook_2010}, but also more recent techniques as Generalized Structured Component Analysis \cite{hwang_generalized_2004}, Generalized Regularized Canonical Correlation Analysis \cite{tenenhaus_regularized_2011} and THEME \cite{bry_theme-seer:_2012}, \cite{bry_theme_2015}.\\
Factor-methods and PLS-type ones have been compared in several works \cite{joreskog_recent_1982}. The gist is that the latter encounter less convergence problems than the former with small samples. A second advantage is that, since they express every LV as a linear combination of its OV's, and yield the estimated coefficients of that combination, the values of the LV's are estimated, and can also be forecast on new samples, opening the way to easy cross-validation. Such is not the case of factor-methods, which do have yet the theoretical advantage to be based on a proper statistical distribution-based model of data, contrary to PLS-type methods, thus allowing standard statistical tests, which are not possible with the latter.\\
In many areas, it is of essence to be able to estimate the values of LV's on statistical units, since these values allow to efficiently analyze the disparities of units on a reduced number of dimensions.\\
Therefore, we are interested in estimating these values even in the factor-model context. In this work, we adapt the EM algorithm to the SEM estimation problem, in order to get estimates of the factor values. 
The paper is organized as follows. Section 2 formally introduces the equations of the SEM we deal with. Section 3 applies the EM algorithm to the SEM and derives the estimation formulas. Section 4 first presents a simulation-based study of the performance of the method, with comparison to more classical methods, and then an application to environmental data.

\section{The model}
\subsection{Notations}
\subsubsection{Data notations}
The data consists in blocks of OV's describing the same $n$ units.
We consider the following data-matrices and notations:\\
$Y= \lbrace y_i^j \rbrace $; $i \in  \llbracket 1,n \rrbracket$, $j \in  \llbracket 1,q_Y \rrbracket$ is the $n \times q_Y$ matrix coding the dependent block of OV's $y^1,...,y^{q_Y}$, identified with its column-vectors.\\
$X^m= \lbrace x_i^{j,m}  \rbrace $; $i \in  \llbracket 1,n \rrbracket$, $j \in  \llbracket 1,q_m \rrbracket$, $m \in  \llbracket 1,p \rrbracket$ is the $n \times q_m$ matrix coding the $m^{ieth}$-explanatory block of OV's $x^{1,m},...,x^{q_m,m}$. Value of variable $x^{j,m}$ for unit $i$ is denoted $x_i^{j,m}$. Variable-blocks will be referred to through the corresponding matrix.\\
$T$ (resp. $T^1,...,T^p$) refers to a $n \times r_T$ (resp. $n \times r_1$, ..., $n \times r_p$) matrices of covariates. \\
We assume that:
\begin{itemize}
\item The units, hence the rows of matrices $Y, X^1,..., X^p$ are independent multivariate normal vectors.
\end{itemize} 
\subsubsection{Model notations}
\begin{figure}[H]
\centering
\includegraphics[width = 1\textwidth]{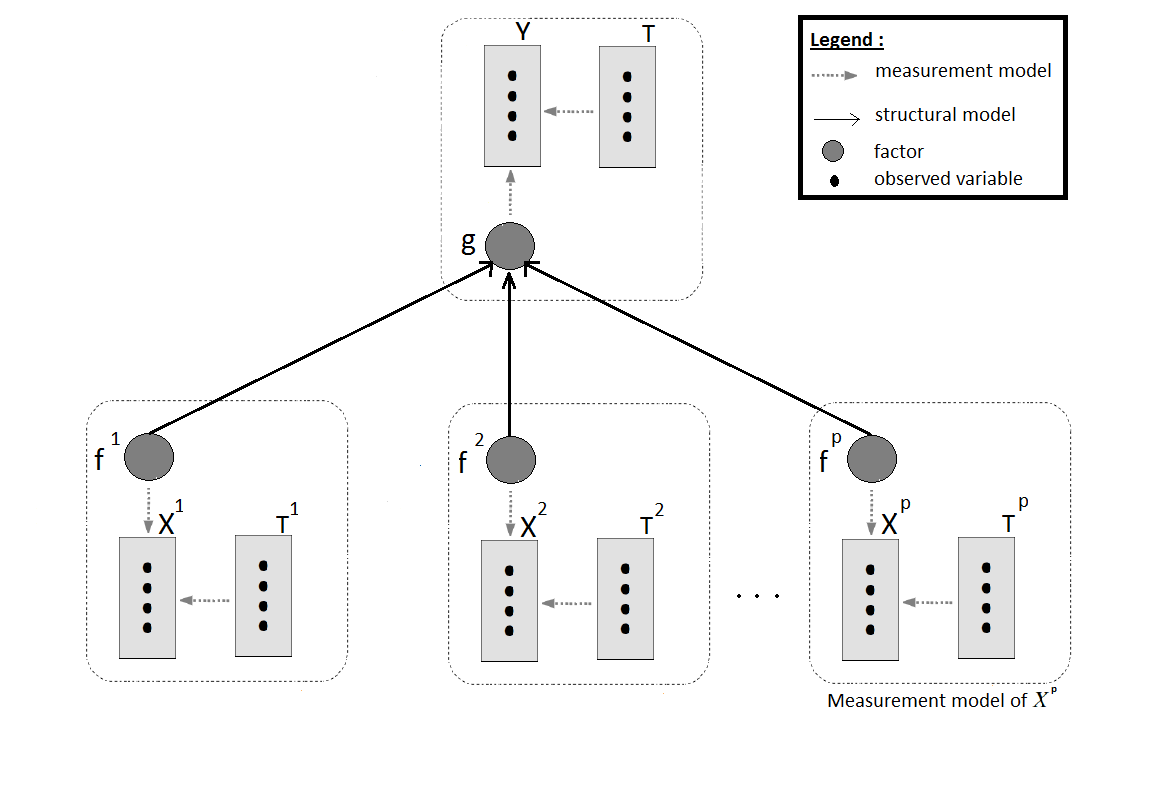}
\caption{Model with a dependent block and $p$ explanatory blocks}
\label{schema_p}
\end{figure}

 The SEM we handle here is a restricted one, in that it contains only one structural equation, relating a dependent factor $g$, underlying a block $Y$ of OV's, to $p$ explanatory factors $f^1,..., f^p$ respectively underlying blocks $X^1,..., X^p$ of OV's (cf. fig.\ref{schema_p}). The main assumptions of this model are the following:
\begin{itemize}
\item Factors $f^1,..., f^p$ are standard normal, i.e. $ \forall m \in \llbracket 1,p \rrbracket, {\mathbb E}(f^m) = 0$ and ${\mathbb V}(f^m)=I_n$.
\item In each block (e.g. $X^p$), the OV's (e.g. $x_j^m, j \in \llbracket 1,q_p \rrbracket$) depend linearly on the block's factor (e.g. $f^m$) and a block of extra-covariates (e.g. $T^m$), conditional on which they are independent.
\item Factor $g$ is normal with zero-mean, and its expectation conditional on $f^1,..., f^p$ is a linear combination of them.
\end{itemize} 
The SEM consists of $p+1$ measurement equations and one structural equation. It is graphed on (cf. fig.\ref{schema_p}).

\subsection{Measurement equations}
As formerly mentioned, each measurement equations relates the variables in a block $X^m$ (respectively $Y$) to the block's factor  $f^m$ (resp. $g$). This link may also involve covariates  $T^m$ (resp. $T$): each OV is expressed as a linear combination of the factor, the covariates and some noise. Hence the $p+1$ measurement equations:
\begin{equation*}
\left\{ \begin{array}{ccl}
Y &=& TD + g b' + \varepsilon^{Y}\\
\forall m \in \llbracket 1,p \rrbracket, 
X^{m} &=& T^m D^m + f^m {a^m}' + \varepsilon^m
\end{array}
\right.
\end{equation*}
where $D$ (resp. $D^m$) is a $r_T \times q_Y$ (resp. $r_m \times q_m$) parameter matrix, $b$ (resp. $a^m$) a $1 \times q_Y$ (resp. $1 \times q_m$) parameter matrix, and $\varepsilon^{Y}$ (resp. $\varepsilon^m$) an $n \times q_Y$ (resp. $n \times q_m$) measurement-error matrix.\\
We impose that the first column of $T$ as well as of each $T^m$ matrix is equal to constant vector having all elements equal to 1. Thus, the first row of $D$ and of each $D^m$ contains mean-parameters.\\
As far as distributions are concerned, we assume that
\begin{itemize}
\item $\varepsilon_{i}^{Y}$ $\sim$ ${\cal N}$(0, $\psi_{_{Y}}$) , where $\psi_{_{Y}} = diag(\sigma^2_{Y,j})_{j \in \llbracket 1,q_{_Y} \rrbracket}$.
\item $\forall m \in \llbracket 1,p \rrbracket$: $\varepsilon_{i}^m$ $\sim$ ${\cal N}$(0, $\psi_m$), where $\psi_m = diag(\sigma^2_{m,j})_{j \in \llbracket 1,q_m\rrbracket}$ and that
\item $\varepsilon^{Y}$ and $\varepsilon^m$, $\forall m \in \llbracket 1,p \rrbracket$ are independent. 
\end{itemize}

As to the factors, we assume that:
\begin{itemize}
\item $\forall m \in \llbracket 1,p \rrbracket$: $f^m$ $\sim$ ${\cal N}(0, Id_n)$ with $f^1,...,f^m$ independent. 
\end{itemize}

\subsection{Structural equations}
The structural equation we consider relates dependent factor $g$ to explanatory factors $f^1,...,f^p$ (cf. fig.\ref{schema_p}) through a linear model:
\begin{equation*}
\left\{ \begin{array}{ccl}
g &=& f^1 c^1+ \cdots +f^{^{p}} c^{^{p}} + \varepsilon^{g}
\end{array}
\right.
\end{equation*}
where $\forall m \in \llbracket 1,p \rrbracket$, $c^m$ is a scalar parameter, and $\varepsilon^{g} \in \mathbb{R}^n$ is a disturbance vector.\\
We assume that
\begin{itemize}
\item $\varepsilon^{g}$ $\sim$ ${\cal N}(0, 1)$
\item $\varepsilon^{g}$ is independent of $\varepsilon^{Y}$ and $\varepsilon^{m}$, $\forall m \in \llbracket 1,p \rrbracket$.
\end{itemize}
N.B. The unit-variance of disturbance $\varepsilon^{g}$ serves an identification purpose.\\
Hence we have the overall model:
\begin{equation}
\left\{ \begin{array}{ccl}
Y &=& TD + g b' + \varepsilon^{Y}\\
\forall m \in \llbracket 1,p \rrbracket, 
X^{m} &=& T^m D^m + f^m {a^m}' + \varepsilon^m \\
g &=& f^1 c^1+ \cdots +f^{^{p}} c^{^{p}} + \varepsilon^{g}
\end{array}
\right.
\label{model_p}
\end{equation}
where the set of parameters is $\theta=\lbrace D, D^1,..., D^p,  b, {a^1},..., {a^p}, {c^1}, {c^2}, \psi_Y, \psi_{1},..., \psi_{p} \rbrace$ such as $\theta$ is $K$-dimensional.\\
Thus, when all $\psi$ matrices are diagonal, we have: 
\begin{align}
\begin{split}
K 
&= 2 + q_Y(r_T+2) + \sum \limits_{\substack{m=1}}^{p} q_m(r_{m}+2)\\
\end{split}
\label{K}
\end{align}

\subsection{A simplified model}
But in order to avoid heavy formulas in the development of the algorithm, we shall use in the sequel, with no loss of generality, a simplified model involving $p=2$ explanatory blocks $X^1$ and $X^2$. The corresponding equation set, for a given unit $i$, reads:
\begin{equation}
\left\{ \begin{array}{ccl}
y_i' &=& {t_i}'D+ g_i b' + {{\varepsilon_i}^{y}}' \\
{x_i^1}'&=& {t_i^1}'D^1 + f_i^1 {a^1}' + {{\varepsilon_i}^{1}}' \\ 
{x_i^2}'&=& {t_i^2}'D^2  + f_i^2 {a^2}' + {{\varepsilon_i}^{2}}' \\ 
g_i &=& f_i^1 c^1 + f_i^2 c^2 + {{\varepsilon_i}^{g}}
\end{array}
\right.
\label{model_p2}
\end{equation}
Such as, $\theta=\lbrace D, D^1, D^2,  b, {a^1}, {a^2}, {c^1}, {c^2},\sigma^2_Y, \sigma^2_1, \sigma^2_2 \rbrace$. Thus, in this case (cf. \eqref{K}), the dimension of $\theta$ is:
\begin{align*}
\begin{split}
K 
&= 5 + q_Y(r_T+1) + \sum \limits_{\substack{m=1}}^{2} q_m(r_{m}+1)\\
\end{split}
\end{align*}

\section{Estimation using the EM algorithm}
In this work, likelihood maximization is carried out via an iterative EM algorithm (\cite{dempster_maximum_1977}, section 4.7). Each iteration of the algorithm involves an Expectation (E)-step followed by a Maximization (M)-step. \cite{dempster_maximum_1977} prove that the EM algorithm yields maximum likelihood estimates. Moreover, they proved that even if the starting point is one where the likelihood is not convex, if an instance of the algorithm converges, it will converge to a (local) maximum of the likelihood.
Another major advantage of the EM algorithm is that it can be used to "estimate" missing data. Thus, if we consider LV's as missing data, the EM algorithm will prove a general technique to maximize the likelihood of statistical models with LV's, but also to estimate these LV's. In our SEM framework, LV's correspond to factors. Thus, we will be able to estimate the factors at unit-level.
We shall present the algorithm on the simplified model (cf. section 2.4) with no loss of generality.

\subsection{The EM algorithm}
Let $z=(y, {x^1}, {x^2})$ be the OV's and $h=(g, f^1, f^2)$ the LV's. The EM algorithm is based on the log-likelihood associated with the complete data $(z,h)$. 
 
\subsubsection{The complete log-likelihood function}
Let $p(z,h;\theta)$ denote the probability density of the complete data. The corresponding log-likelihood function is:
\begin{align*}
\begin{split}
{\cal L}(\theta;z,h) = - \frac{1}{2} 
\sum \limits_{\substack{i=1}}^{n}
&\lbrace
ln|\psi_Y| + ln|\psi_{1}| + ln|\psi_{2}| \\
&+ (y_i - D't_i - g_ib)'\psi_Y^{-1}(y_i - D't_i - g_ib) \\
&+ (x_i^1 - {D^1}'t_i^1 - f_i^1a^1)'\psi_{1}^{-1}(x_i^1 - {D^1}'t_i^1 - f_i^1a^1) \\
&+ (x_i^2 - {D^2}'t_i^2 - f_i^2a^2)'\psi_{2}^{-1}(x_i^2 - {D^2}'t_i^2 - f_i^2a^2) \\
&+ (g_i - c^1\ f_i^1 - c^2\ f_i^2)^2 + (f_i^1)^2 + (f_i^2)^2 \rbrace + \lambda
\\
\end{split}
\end{align*}
Where $\theta$ is the $K$-dimensional set of model parameters and $\lambda$ a constant. However, because of the simplification made in the section 2.4, in our case  $\theta=\lbrace D, D^1, D^2,  b, {a^1}, {a^2}, {c^1}, {c^2},\sigma^2_Y, \sigma^2_1, \sigma^2_2 \rbrace$. Indeed, $\psi_Y=\sigma^2_Y Id_{q_Y}$, $\psi_{1}=\sigma^2_1 Id_{q_1}$ and $\psi_{2}=\sigma^2_2 Id_{q_2}$.

\subsubsection{Estimation in the SEM}
To maximize this function, in the framework of EM algorithm, we have to solve:
\begin{equation}
\mathbb{E}_z^h[\frac{\partial}{\partial\theta}{\cal L}(\theta;z,h)] = 0
\label{equation_EM}
\end{equation}
\cite{foulley_algorithme_2002}

This demands that we know the derivatives of the log-likelihood function and the distribution $p_{z_i}^{h_i}$ of $h_i$ conditional on $z_i$ for each observation $i \in \llbracket 1,n \rrbracket$. Let us introduce the following notation:\\
\begin{center}
$p_{z_i}^{h_i} = {\cal N}(M_i = \begin{pmatrix} {m_1}_i \\ {m_2}_i \\ {m_3}_i \end{pmatrix} , \Sigma = \begin{pmatrix} \sigma_{{11}} & \sigma_{{12}} & \sigma_{{13}}\\ \sigma_{{21}} & \sigma_{{22}} & \sigma_{{23}} \\ \sigma_{{31}} & \sigma_{{32}} & \sigma_{{33}} \end{pmatrix})$
\end{center}

\begin{center}
\begin{tabular}{ll} 
$\widetilde{g_i} = \mathbb{E}_{z_i}^{h_i} [ g_i] = {m_1}_i$ ; & $\widetilde{{\cal \gamma}_i} = \mathbb{E}_{z_i}^{h_i}[g_i^2] = (\mathbb{E}_{z_i}^{h_i}[g_i])^2 + \mathbb{V}_{z_i}^{h_i}[g_i] = {m_1}_i^2 + \sigma_{{11}}$\\
$\widetilde{f_i^1} = \mathbb{E}_{z_i}^{h_i} [ f_i^1] = {m_2}_i$ ; & $ \widetilde{{\cal \phi}_i^1} = \mathbb{E}_{z_i}^{h_i}[(f_i^1)^2] = (\mathbb{E}_{z_i}^{h_i}[f_i^1])^2 + \mathbb{V}_{z_i}^{h_i}[f_i^1] = {m_2}_i^2 + \sigma_{{22}} $\\
$\widetilde{f_i^2} = \mathbb{E}_{z_i}^{h_i} [ f_i^2] = {m_3}_i$ ; & $ \widetilde{{\cal \phi}_i^2}  = \mathbb{E}_{z_i}^{h_i}[(f_i^2)^2] = (\mathbb{E}_{z_i}^{h_i}[f_i^2])^2 + \mathbb{V}_{z_i}^{h_i}[f_i^2] = {m_3}_i^2 + \sigma_{{33}}$\\
\end{tabular}
\end{center}
For all $\widetilde{\xi} \in \lbrace \widetilde{g}, \widetilde{f^1}, \widetilde{f^2}, \widetilde{{\cal \gamma}}, \widetilde{{\cal \phi}^1}, \widetilde{{\cal \phi}^2} \rbrace$, we denote $\widetilde{\xi} = {(\widetilde{\xi}_i)}_{i=1,...,n} \in \mathbb{R}^n$.\\
The parameters of the gaussian distribution $p_{z_i}^{h_i}$ are explicit and have the following form:
\begin{center}
$ M_i = \Sigma^*_2 {\Sigma^*_3}^{-1} \mu^*$ and $\Sigma = \Sigma^*_1 - \Sigma^*_2 {\Sigma^*_3}^{-1}{\Sigma^*_2}'$ where:
\end{center}
$\Sigma^*_1 = \begin{pmatrix} (c^1)^2+ (c^2)^2+1 & c^1 & c^2 \\  
                                 c^1 & 1 & 0 \\
                                 c^2 & 0 & 1
                                 \end{pmatrix}$\\
$\Sigma^*_2 = \begin{pmatrix} ((c^1)^2+ (c^2)^2 + 1)b'  &  c^1{a^1}' & c^2{a^2}' \\  
                        c^1b' &  {a^1}' & 0_{(1,q_2)}\\
                        c^2b' & 0_{(1,q_1)} & {a^2}' \end{pmatrix}$\\
$\Sigma^*_3 = \begin{pmatrix} ((c^1)^2+ (c^2)^2 + 1)bb' + \Psi_{Y} & c^1b{a^1}' & c^2b{a^2}' \\
                       c^1a^1b' & a^1{a^1}' + \Psi_{1} &  0_{(q_1,q_2)}\\
                      c^2a^2b' & 0_{(q_2,q_1)} & a^2{a^2}' + \Psi_{2}
                       \end{pmatrix}$\\                       
$\mu^*_i = \begin{pmatrix} y_i - D't_i  \\
                       x_i^1 - {D^1}'t_i^1 \\
                       x_i^2 - {D^2}'t_i^2 \end{pmatrix}$\\
These results are demonstrated in AppendixB. Expressions of the first-order derivatives of ${\cal L}$ with respect to $\theta$ are also established in AppendixB and written in the following forms with $m \in \lbrace 1,2 \rbrace$:

\begin{equation}
\left\{ \begin{array}{ccl}
\frac{\partial}{\partial D'}{\cal L}(z,h) &=& \sum \limits_{\substack{i=1}}^{n}
 \psi_Y^{-1}(y_i - D't_i - g_ib) {t_i}' \\
 \frac{\partial}{\partial {D^{m}}'}{\cal L}(z,h) &=& \sum \limits_{\substack{i=1}}^{n}
\psi_{m}^{-1}(x_i^m - {D^m}'t_i^m - f_i^ma^m){t_i^m}'\\
\frac{\partial}{\partial b}{\cal L}(z,h) &=& \sum \limits_{\substack{i=1}}^{n}
g_i\psi_Y^{-1}(y_i - D't_i - g_ib)\\
\frac{\partial}{\partial a^m}{\cal L}(z,h) &=& \sum \limits_{\substack{i=1}}^{n}
f_i^m\psi_{m}^{-1}(x_i^m - {D^m}'t_i^m - f_i^ma^m)\\
\frac{\partial}{\partial c^m}{\cal L}(z,h) &=& \sum \limits_{\substack{i=1}}^{n}
f_i^m(g_i- c^2\ f_i^2 - c^1\ f_i^1)\\
\frac{\partial}{\partial\sigma_Y^2}{\cal L}(z,h) &=& n\ q_{_Y}\ \sigma_Y^{-2}-\sigma_Y^{-4} 
\sum \limits_{\substack{i=1}}^{n}
|| y_i - D't_i - g_ib ||^2\\
\frac{\partial}{\partial\sigma_{m}^2}{\cal L}(z,h) &=& n\ q_{_{m}}\ \sigma_{m}^{-2}-\sigma_{m}^{-4} 
\sum \limits_{\substack{i=1}}^{n}
|| x_i^m - {D^m}'t_i^m - f_i^ma^m ||^2
\end{array}
\right.
\label{syst_derrivees1}
\end{equation}

So, here formula \eqref{equation_EM} develops into:

\begin{equation}
\hspace{-1.8cm}
\left\{ \begin{array}{ccl}
\sum \limits_{\substack{i=1}}^{n}
(y_i-D't_i - \widetilde{g_i}b) {t_i}' &=& 0\\
\sum \limits_{\substack{i=1}}^{n}
(x_i^m - {D^m}'t_i^m  -\widetilde{f_i^m}a^m){t_i^m}' &=& 0\\
\sum \limits_{\substack{i=1}}^{n}
\widetilde{g_i}y_i-\widetilde{g_i}D't_i - \widetilde{{\cal \gamma}_i}b &=& 0\\
\sum \limits_{\substack{i=1}}^{n}
\widetilde{f_i^m}x_i^m-\widetilde{f_i^m}{D^m}'t_i^m - \widetilde{{\cal \phi}_i^m}a^m &=& 0\\
\sum \limits_{\substack{i=1}}^{n}
\sigma_{12} + \widetilde{f_i^1}\widetilde{g_i} - c^2\ \sigma_{23} - c^2\ \widetilde{f_i^1}\widetilde{f_i^2} - \widetilde{{\cal \phi}_i^1}c^1 &=& 0\\
\sum \limits_{\substack{i=1}}^{n}
\sigma_{31} + \widetilde{f_i^2}\widetilde{g_i} - c^2\ \widetilde{\phi_i^2} - c^1 \sigma_{32}  - c^1\widetilde{f_i^1}\widetilde{f_i^2}  &=& 0\\
nq_{_Y}\sigma_Y^{-2}-\sigma_Y^{-4} 
\sum \limits_{\substack{i=1}}^{n}
|| y_i-D't_i||^2 + ||b||^2\widetilde{{\cal \gamma}_i} -2(y_i-D't_i)'\widetilde{g_i}b &=& 0\\
nq_m\sigma_m^{-2}-\sigma_m^{-4} 
\sum \limits_{\substack{i=1}}^{n}
|| x_i^m-{D^m}'t_i^m||^2 + ||a^m||^2\widetilde{{\cal \phi}_i^m} -2(x_i^m-{D^m}'t_i^m)'\widetilde{f_i^m}a^m &=& 0
\end{array}
\right.
\label{syst_derrivees2}
\end{equation}

System of equations \eqref{syst_derrivees2} is easy to solve and the obtained solutions will be given in the next section.

\subsubsection{Results}
The explicit solution of the system \eqref{syst_derrivees2} (and also of \eqref{equation_EM}) is the following:
\vspace{0.2cm}

\fbox{
\begin{minipage}{380pt}
\begin{align}
\begin{split}
& \hat{b} 
= \frac{(\overline{\mathstrut \widetilde{ g} y } - \overline{\mathstrut yt'})(\overline{tt'})^{-1} \  \overline{\mathstrut \widetilde{ g}t}}{\overline{\mathstrut \widetilde{ {\cal \gamma}}} - \overline{\mathstrut \widetilde{ g}t'}(\overline{tt'})^{-1}\overline{\mathstrut \widetilde{ g}t}}\\
& \widehat{a^m} 
= \frac{\overline{\mathstrut \widetilde{ f^m} x^m } - \overline{\mathstrut x^m{t^m}'}( \overline{t^m{t^m}'})^{-1} \overline{\mathstrut \widetilde{ f^m}t^m}}{\overline{\mathstrut \widetilde{ {\cal \phi}^m}} - \overline{\mathstrut \widetilde{ f^m}{t^m}'}( \overline{t^m{t^m}'})^{-1}\overline{\mathstrut \widetilde{ f^m}t^m}}\\
& \widehat{c^1} 
= \frac{(\sigma_{12} + \overline{\widetilde{ f^1}\widetilde{ g}}){\overline{\widetilde{\phi}^2}} - ( \sigma_{13} + \overline{\widetilde{ f^2}\widetilde{ g}} )( \sigma_{23} + \overline{\widetilde{ f^1}\widetilde{ f^2}} )
}{\overline{\widetilde{\phi}^1} \overline{\widetilde{\phi}^2} - ( \sigma_{23} + \overline{\widetilde{ f^1}\widetilde{ f^2}} ) ^2}\\
& \widehat{c^2} 
= \frac{(\sigma_{13} + \overline{\widetilde{ f^2}\widetilde{ g}}){\overline{\widetilde{\phi}^1}} - ( \sigma_{12} + \overline{\widetilde{ f^1}\widetilde{ g}} )( \sigma_{23} + \overline{\widetilde{ f^1}\widetilde{ f^2}} )
}{\overline{\widetilde{\phi}^1} \overline{\widetilde{\phi}^2} - ( \sigma_{23} + \overline{\widetilde{ f^1}\widetilde{ f^2}} ) ^2}\\
&\widehat{D'} 
=(\overline{\mathstrut yt'} - \widehat{b} \ \overline{\mathstrut\widetilde{g}t'})(\overline{tt'})^{-1}\\
&\widehat{{D^m}'} 
=(\overline{\mathstrut x^m{t^m}'} - \widehat{a^m} \overline{\mathstrut\widetilde{f^m} {t^m}'})( \overline{t^m{t^m}'})^{-1}\\
&\widehat{\sigma_Y^2} 
= \frac{1}{nq_{_Y}} \sum \limits_{\substack{i=1}}^{n}
\lbrace
||y_i- \widehat{D'}t_i||^2+||\hat{b}||^2\widetilde{{\cal \gamma}_i}-2(y_i-\widehat{D'}t_i)'\hat{b}\widetilde{ g_i}
\rbrace \\
&\widehat{\sigma_m^2} 
= \frac{1}{nq_m} \sum \limits_{\substack{i=1}}^{n}
\lbrace
||x_i^m-\widehat{{D^m}'}t_i^m||^2+||\widehat{a^m}||^2\widetilde{{\cal \phi}_i^m}-2(x_i^m-\widehat{{D^m}'}t_i^m)'\widehat{a^m}\widetilde{ f_i^m}
\rbrace \\
\end{split}
\label{formules_solution}
\end{align}
\end{minipage} }

\subsubsection{The algorithm}
 To estimate parameters in $\theta$, we propose the following EM-algorithm. We denote $[t]$ the $t^{ieth}$-iteration of the algorithm.
\vspace*{1cm}

\begin{enumerate}
\item{Initialization \footnote{In the initialization step, $\forall m \in \lbrace 1,p \rbrace$ we propose to obtain ${D^m}^{[0]}$ by multiple linear regression between $X^m$ and $T^m$. Then, to initialize the others, we compute each approximated factor $\widetilde{f^m}^{[0]}$ and $\widetilde{g}^{[0]}$ as first principal component of a PCA  of $X^m - T^m {D^m}^{[0]}$ and  $Y - T {D}^{[0]}$. Thus, we initialize $a^m$, $\sigma^2_{m}$ (resp. $b$, $\sigma^2_y$) by multiple linear regression between $X^m - T^m {D^m}^{[0]}$ and $\widetilde{f^m}^{[0]}$ (resp. between $Y - T {D}^{[0]}$ and $\widetilde{g}^{[0]}$). Finally, each ${c^m}^{[0]}$ can be obtained by multiple linear regression between $\widetilde{g}^{[0]}$ and $\sum \limits_{\substack{m=1}}^{p} \widetilde{f^m}^{[0]}$. In practice we use the functions \textit{lm()} and \textit{PCA()} derived from the package \textit{FactoMineR} \cite{husson_factominer:_2008}.} = choice of the initial parameter values $\theta^{[0]}$.}\\ 
\item{Current iteration $t \geq 1$, until stopping condition is met:
\begin{enumerate}
\item \textbf{E-step}: with $\theta^{[t-1]}$,
\begin{enumerate}
\item Calculate explicitly distribution $p_{z_i}^{h_i}$ for each $i \in \llbracket 1,n \rrbracket$.
\item Estimate the factor-values  $\widetilde{g}^{[t]}$, $\widetilde{f^m}^{[t]}$, $m \in \lbrace 1, 2 \rbrace$.
\item Calculate $\widetilde{{\cal \gamma}}^{[t]}$ and $\widetilde{{\cal \phi}^m}^{[t]}$, $m \in \lbrace 1, 2 \rbrace$.
\end{enumerate}
\item \textbf{M-step}:
\begin{enumerate}
\item Update  $\theta$ to $\theta^{[t]}$ by injecting $\widetilde{g}^{[t]}$, $\widetilde{{\cal \gamma}}^{[t]}$ and $\widetilde{f^m}^{[t]}$,  $\widetilde{{\cal \phi}^m}^{[t]}$, $m \in \lbrace 1,2 \rbrace$ into the formulas in \eqref{formules_solution}.
\end{enumerate}
\end{enumerate}
}
\item{ We used the following stopping condition with the smallest $\epsilon$ possible:}
\end{enumerate}

\begin{equation}
\sum \limits_{\substack{k=1}}^{K}  \dfrac{\vert {\theta^{*}}^{[t+1]}[k] - {\theta^{*}}^{[t]}[k] \vert}{{\theta^{*}}^{[t+1]}[k]} < \epsilon
\end{equation}
where $\theta^{*}$ is the $K$-dimensional vector containing the scalar values in all parameters in $\theta$.

\section{Numerical results on simulated data}
\subsection{Data generation}
We consider $n=400$ units and 
$q_Y=q_1=q_2= 40$. Therefore, the $120$ OV's $Y,X^1,X^2$ are simulated so as to be structured respectively around three factors $g,f^1,f^2$. Factors $f^1$ and $f^2$ are explanatory of $g$. Besides, we consider $r_T=r_{1}=r_{2}=2$ i.e $2$ covariates are simulated for each covariate matrix $T$, $T^1$ and $T^2$. The data is simulated as follows.
\begin{enumerate}
\item{Choice of $\theta$:}
\begin{enumerate}
\item $D=D^1=D^2$ a) matrices filled in row-wise with the ordered integer sequence ranging from $1$ to $80$ (indeed: $r_T*q_Y=r_{1}*q_{1}=r_{2}*q_{2}=2*40$).
\item $b=a^1=a^2=$ ordered integer sequence ranging from $1$ to $40$.
\item $c^1=c^2=1$
\item $\sigma^2_{Y}=\sigma^2_{1}=\sigma^2_{2}=1$
\end{enumerate}
\item{Simulation of factors $g,f^1,f^2$}
\begin{enumerate}
\item Simulate vectors $f^1$ and $f^2$ of $n=400$ normally distributed random numbers with mean 0 and variance 1 (abbreviated $\forall m ,\in \lbrace 1,2 \rbrace$ $f^m$ $\sim$ ${\cal N}(0, Id_{400})$).
\item We simulate ${\varepsilon}^{g}$ according to distribution ${\varepsilon}^{g}$ $\sim$ ${\cal N}(0, Id_{400})$.
\item We then calculate $g$ as $g = f^1 c^1 + f^2 c^2 + {{\varepsilon}^{g}}$
\end{enumerate}
\item{Simulation of noises $\varepsilon^{Y}$, $\varepsilon^{1}$, $\varepsilon^{2}$}
Each element of matrix $\varepsilon^{Y}$, (respectively $\varepsilon^{1}$, $\varepsilon^{2}$) is simulated independently from distribution ${\cal N}(0, \sigma^2_{Y} = 1)$ (respectively, $\sigma^2_1 = 1$, $\sigma^2_2 = 1$). 
\item{Simulation of covariate matrices $T$, $T^1$, $T^2$}\\
Each element of matrices $T$, $T^1$, $T^2$ is simulated according to the standard normal distribution. 
\item{Construction of $Y$, $X^1$, $X^2$}
$Y$, $X^1$, $X^2$ are eventually calculated through formulas in the model \eqref{model_p}.
\end{enumerate}
This simulation scheme was performed 100 times, each time yielding a set of simulated data matrices $(Y,X^1,X^2)$.  Then for each simulated data, we ran an estimation routine with a threshold value $\varepsilon=10^{-2}$, yielding the average results presented in section 4.2. Thus from $400*120=48000$  scalar elements of data, we will estimate $3*n=1200$ scalar elements of factors plus $K=5+3*40(2+1)=365$ scalar parameters, i.e: $1565$ scalars.
\subsection{Results}
\begin{figure}[H]
\centering
\includegraphics[width = 1\textwidth]{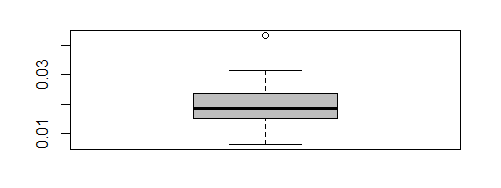}
\caption{Box plot of the average absolute relative deviations between simulated parameters and their estimations.}
\label{boxplot_parameters}
\end{figure}

\begin{figure}[H]
\centering
\includegraphics[width = 1\textwidth]{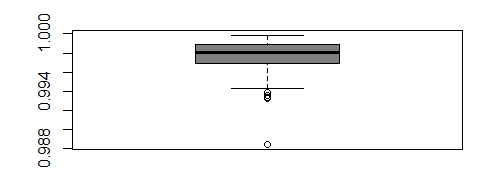}
\caption{Box plot of the correlations between the simulated factors and their estimations}
\label{boxplot_factors}
\end{figure}

Convergence was observed in almost all cases in less than five iterations.
We assess the quality of the estimations as follows.
\begin{itemize}
\item On the one hand, we calculate the absolute relative deviation between each simulated scalar parameter in $\theta^*$ and its estimation, and then average these deviations over the $100$ simulations. We then produce a box-plot of the average absolute relative deviations (cf. fig. \ref{boxplot_parameters}). This makes the interpretation easier, since we only need to look at the box-plot's values and check that they are positive (because of the absolute value) and close to $0$.
\item On the other hand, to assess the quality of the factor estimations, we compute the $300$ values of square correlations between the simulated concatenated factors $(g,f^1,f^2)$ (respectively) and the corresponding estimations ($(\widetilde{g},\widetilde{f^1},\widetilde{f^2})$). Once again, we produce a box-plot of these correlations (cf. fig. \ref{boxplot_factors}) and check that it indicates values close to $1$.
\end{itemize}
Figures \ref{boxplot_parameters} and \ref{boxplot_factors} show clearly that the estimations are very close to the actual quantities. Indeed, on figure \ref{boxplot_parameters}, the median of average absolute relative deviations is $0.018$, the first quartile is $0.015$ and the third quartile is $0.023$. On figure \ref{boxplot_factors}, the median of square correlations is $0.998$, the first quartile is $0.997$ and the third quartile is $0.999$. So, factor $g$ (respectively $f^1$ and $f^2$) turn out to be drawn towards the principal direction underlying the bundles made up by observed variables $Y$ (respectively $X^1$ and $X^2$). Now, we may legitimately wonder how the quality of estimations could be affected by the number of observations and the number of OV's in each block. In the following section we give a sensitivity analysis performed to investigate this.


\subsection{Sensitivity analysis of estimations}
We performed a sensibility analysis on the simulated data presented in section 4.1. The purpose was to study the influence of the block-sizes ($n$, $q_Y$, $q_1$, $q_2$) on the quality of estimation, both of the parameters and the  factors.  To simplify the analysis, we imposed $q_Y=q_1=q_2=q$ and varied $n$ and $q$ separately, i.e. studied the cases $n = 50, 100, 200, 400$ with $q = 40$ and $q = 5, 10, 20, 40$ with $n = 400$. Each case was simulated $100$ times. Therefore, we simulated $800$ data-sets.

\subsubsection{Sensitivity with respect to the number $n$ of observations}
\begin{figure}[H]
\centering
\includegraphics[width = .5\textwidth]{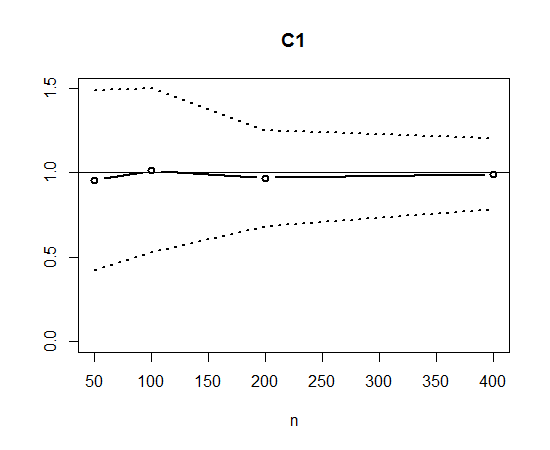}
\caption{Average estimate of $c^1$ and $95\%$ confidence interval as a function of $n$.}
\label{c1_n}
\end{figure}

\begin{figure}[H]
\centering
\includegraphics[width = .5\textwidth]{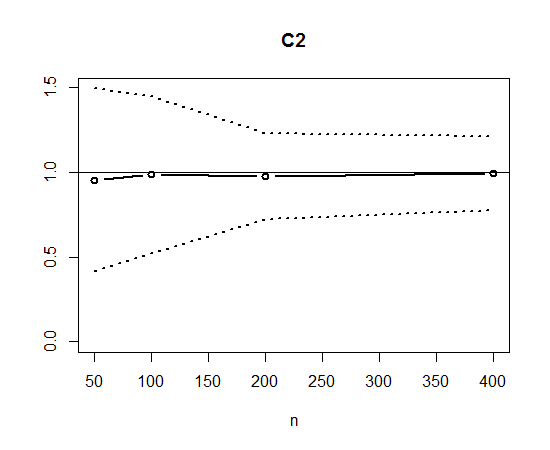}
\caption{Average estimate of $c^2$ and $95\%$ confidence interval as a function of $n$.}
\label{c2_n}
\end{figure}

\begin{figure}[H]
\centering
\includegraphics[width = .5\textwidth]{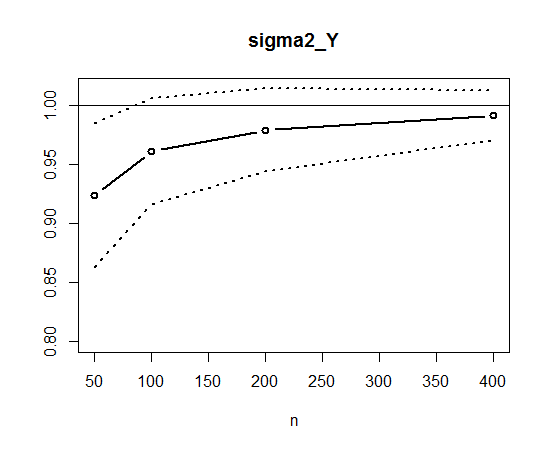}
\caption{Average estimate of $\sigma^2_Y$ and $95\%$ confidence interval as a function of $n$.}
\label{sigma2_n}
\end{figure}

In this section, we study the evolution with $n$ of the average estimation of structural coefficients $c^1$ and $c^2$ and parameter $\sigma^2_Y$ with respect to their actual values, all equal to $1$, and that of the correlations between factors and their estimates. The number of OV's is fixed to $q=40$ in each block. Figures \ref{c1_n}, \ref{c2_n} and \ref{sigma2_n} graph these evolutions (average value of estimate in plain line), including a $95\%$ confidence-interval about each average estimate (dotted line). These figures show that the biases and the standard deviations are, as expected, more important for little values of $n$, but also that the quality of estimation is already quite good for $n=50$. 


\begin{figure}[H]
\centering
\includegraphics[width = 1\textwidth]{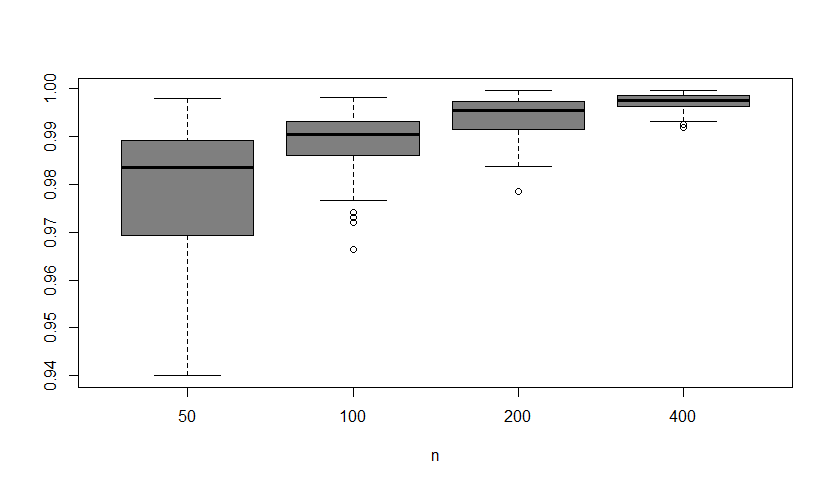}
\caption{Box plots of the correlations between simulated factors and their estimations for various values of $n$.}
\label{boxplot_factors_n}
\end{figure}

As for the factors, figure \ref{boxplot_factors_n} shows that their correlations increase and get close to one as $n$ increases, with a dispersion decreasing to 0. However, even for $n=50$, the correlations are  mostly above $0.95$, indicating that the factors are correctly reconstructed.

\subsubsection{Sensitivity with respect to the number $q$ of OV's in each block}
\begin{figure}[H]
\centering
\includegraphics[width = .5\textwidth]{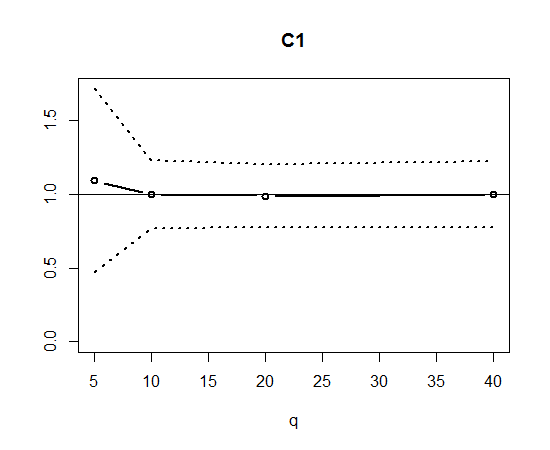}
\caption{Evolution of the average of estimated parameters $c^1$ according to various $n$ values.}
\label{c1_q}
\end{figure}

\begin{figure}[H]
\centering
\includegraphics[width = .5\textwidth]{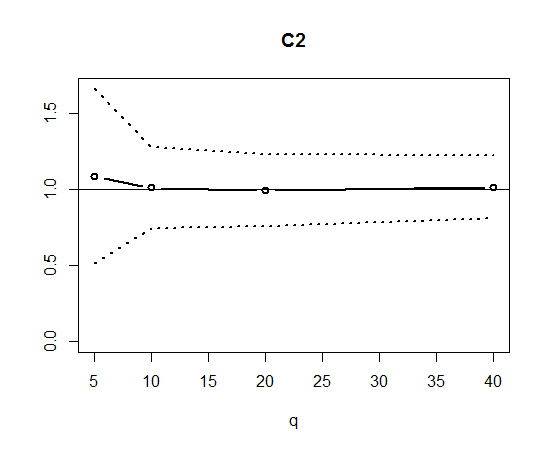}
\caption{Evolution of the average of estimated parameters $c^2$ according to various $n$ values.}
\label{c2_q}
\end{figure}

\begin{figure}[H]
\centering
\includegraphics[width = .5\textwidth]{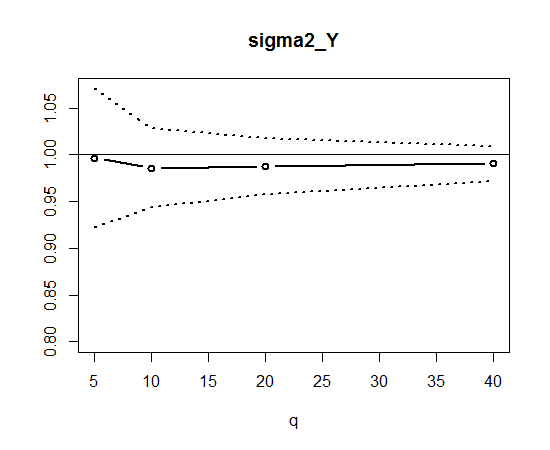}
\caption{Evolution of the average of estimated parameters $\sigma^2_Y$ (at right) according to various $n$ values.}
\label{sigma2_q}
\end{figure}

Likewise, we study the evolution of the average estimates of $c^1$, $c^2$, $\sigma^2_Y$ and the correlation between factors and their estimates for different values of $q$, with  $n$ fixed to $400$. We observe that, unsurprisingly, the biases and the standard deviations decrease as $q$ increases (cf. figures \ref{c1_q}, \ref{c2_q} and \ref{sigma2_q}). We observe that they stabilize even faster with $q$ than with $n$, particularly $\sigma^2_Y$. Indeed, from $q=10$ on, the confidence interval is narrow enough.


\begin{figure}[H]
\centering
\includegraphics[width = 1\textwidth]{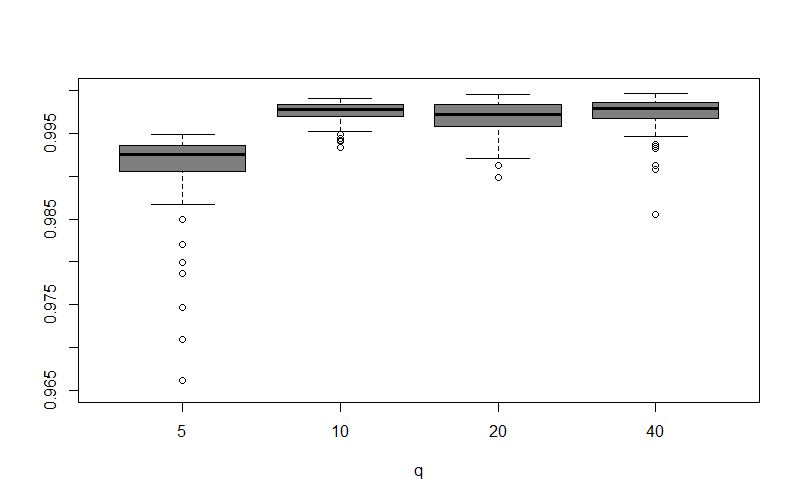}
\caption{Box plots of the correlations between simulated factors and their estimations according to various $n$ values.}
\label{boxplot_factors_q}
\end{figure}

As for the factors, figure \ref{boxplot_factors_q} shows that their correlations are already very close to $1$ for $q=5$, with a very small variance, and keep increasing with $q$. 

To sum things up, the sample size $n$ proved to have more impact on the quality of parameter  estimation and factor reconstruction than the number of OV's. Now, the quality of factor reconstruction remains high for rather small values of $n$ or $q$. We advise to use a minimal sample size of $n=100$ to obtain really stable structural coefficients. Above this threshold, $n$ has but little impact on the biases and standard deviations of estimations.

\section{An application to environmental data}
\subsection{Data presentation}
\begin{figure}[H]
\centering
\includegraphics[width = .75\textwidth]{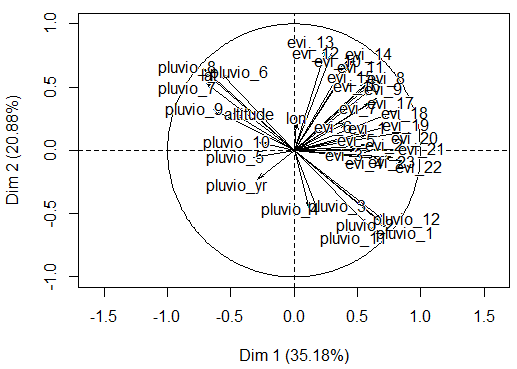}  
\caption{Correlation-scatterplot yielded by the PCA of the $X^1$ and $X^2$ geo-referenced environmental variables.}
\label{ACP}
\end{figure}

We apply our model to the data-set \textit{genus} provided in the R-package SCGLR by \cite{mortier_scglr-r_2014}. Data-set \textit{genus} was built from the CoForChange database. It gives the abundances of 27 common tree genera present in the tropical moist forest of the Congo-Basin, and the measurements of 40 geo-referenced environmental variables, for $n=1000$ inventory plots (observations). Some of the geo-referenced environmental variables describe 16 physical factors pertaining to topography, geology and rainfall description. The remaining variables characterize vegetation through the enhanced vegetation index (EVI) measured on 16 dates.\\ 
In this section, we aim at modeling the tree abundances from the other variables, while reducing the dimension of data. The dependent block of variables $Y$ therefore consists of the $q_Y = 27$ tree species counts divided by the plot-surface. A PCA of the geo-referenced environmental variables and the photosynthetic activity variables confirms that EVI measures are clearly separated from the other variables (cf. Fig. \ref{ACP}). Indeed, Fig. \ref{ACP} shows two variable-bundles with almost orthogonal central directions. This justifies using our model (cf. section 5.2) with $p = 2$ explanatory groups, one of them ($X^1$) gathering $q_1 = 16$ rainfall measures and location variables (longitude, latitude and altitude), and the second one ($X^2$), the $q_2 = 23$ EVI measures. Besides, in view of the importance of the geological substrate on the spatial distribution of tree species in the Congo Basin, showed by \cite{fayolle_geological_2012}, we chose to put nominal variable \textit{geology} in a block $T$. This block therefore contains constant 1 plus all the indicator variables of geology but one, which will be the reference value. Geology having 5 levels, T has thus 5 columns.

\subsection{Model with geologic covariates}
\subsubsection{Model specification}
Here is the model used with the variable-blocks designed in section 5.1.:
\begin{equation}
\left\{ \begin{array}{ccl}
Y &=& TD + g b' + \varepsilon^{Y} \nonumber\\
X^1 &=& \mathds {1}_n {d^1}' + f^1 {a^1}' + \varepsilon^{1} \nonumber\\ 
X^2 &=& \mathds {1}_n {d^2}' + f^{2} {a^2}' + \varepsilon^{2} \nonumber\\
g &=& f^1 c^1 + f^2 c^2 + \varepsilon^{g} \nonumber
\end{array}
\right.
\end{equation}

Where $n=1000$, $q_Y=27$, $q_1=16$, $q_2=23$ and $r_T=5$. The first row of $D$ is a parameter vector that contains the means of the $Y$'s noted $D[1,]$ in Table \ref{tabl_estim_param_db_avec}, and the other rows, the overall effects of the geological substrates with respect to the reference one. Indeed, the next section presents the model's parameter-estimations where, in Table \ref{tabl_estim_param_db_avec}, each row $r$ of $D$ is noted $D[r,]$.

\subsubsection{Results}
With a threshold value $\varepsilon=10^{-3}$, convergence was reached after 58 iterations.
Some parameter-estimations are presented in Tables \ref{tabl_estim_param_db_avec}, \ref{tabl_estim_param_d1a1_avec} and \ref{tabl_estim_param_scal_avec}. For practical reasons, in Appendix D table \ref{tabl_estim_param_d2a2_avec} present the remainder of the parameter-estimations.

\begin{center}
\begin{tabular}{|C{1.5cm}|C{1cm}C{1cm}C{1cm}C{1cm}C{1cm}|C{1cm}||C{3
cm}|}
  \hline
   & \multicolumn{6}{c||}{Parameter-estimations} & \\
   \hline
 Variables & $D[1,]$ & $D[2,]$ & $D[3,]$ & $D[4,]$ & $D[5,]$ & $b'$ & Correlations with $\widetilde{g}$  \\
  \hline
  gen1 & 0.76  & 0.16 & 0.06 & 0.68 & -0.12 & -0.13 & -0.14 \\
  gen2 &  0.54 & -0.28 & -0.03 & -0.03 & -0.28 & 0.47 & 0.58 \\
  gen3 & 0.41  & -0.23 & -0.02 & 0.25 & -0.37 & 0.29 & 0.36\\
  gen4 & 0.12 & 0.14 & 0.03 & 0.52 & 0.30 &  0.09 & 0.15\\
  gen5 & 0.31 & 0.15 & 0.19 & -0.20 & 0.84 & 0.09 & 0.16\\
  gen6 & 0.55 & -0.12 & -0.26 & 0.06 & -0.02 &  0.14 & 0.18\\
  gen7 & 0.46 & 0.06 & -0.04 & -0.37 & 0.43 & 0.14 & 0.18\\
  gen8 & 0.55 & 0.04 & -0.09 & -0.16 & 0.04 & 0.42 & 0.52\\
  gen9 & 0.92 & -0.54 & 0.26 & -0.66 & -0.61 & 0.07 & 0.03\\
  gen10 & 0.68 & 0.40 & 0.20 & 0.37 & 0.06 & -0.32  & -0.39\\
  gen11 & 1.74 & -0.50 & -0.21 & 0 & -0.67 & 0.33 & 0.39\\
  gen12 & 0.87 & 0.14 & 0.73 & -0.51 & -0.21 & 0.24 & 0.26\\
  gen13 & 1.08 & -0.09 & -0.37 & -0.02 & -0.53  & 0.26 & 0.29\\
  gen14 & 0.41 & -0.16 & -0.10 & 0.12 & -0.36 & -0.05 & -0.07\\
  gen15 & 0.51 & 0.01 & -0.11 & 0.27 & -0.18 & 0.29 & 0.37\\
  gen16 & 0.50 & -0.19 & -0.01 & 0.55 & -0.27 & 0.1 & 0.14\\
  gen17 & 0.79& -0.54 & -0.20 & -0.52 & -0.45 & 0.39 & 0.45\\
  gen18 & 0.16 & -0.05 & 0.20 & 0.03 & -0.03 & 0.18 & 0.23\\
  gen19 & 0.34 & 0.06 & 0.41 & -0.11 & 0.38 & 0.23 & 0.31\\
  gen20 & 0.49 & 0.02 & -0.21 & 0.08 & 0.14 & -0.2 & -0.24\\
  gen21 & 0.79 & -0.30 & -0.12 & 0.71 & -0.13 & 0.12 & 0.19\\
  gen22 & 0.32 & -0.07 & -0.07 & 0.38 & -0.11 &  0.23 & 0.3\\
  gen23 & 1.02 & -0.28 & -0.31 & 0 & -0.07 & 0.46 & 0.58\\
  gen24 & 0.80 & -0.23 & -0.08 & 0.22 & -0.47 & 0.57 & 0.7\\
  gen25 & 0.60 & -0.16 & -0.04 & 0.97 & -0.49 & 0.41 & 0.53\\
  gen26 & 0.84 & 0.22 & 0.27 & -0.70 & 0.82 & 0.04 & 0.07\\
  gen27 & 0.27 & 0.41 & 0.69 & -0.24 & 0.56 & 0.08 & 0.11\\
 \hline
\end{tabular}
\captionof{table}{Application to the \textit{genus} data with geologic covariate : estimations of parameters $D'$ and $b'$, and correlations between $\widetilde{g}$ and the variables $Y$}
\label{tabl_estim_param_db_avec}

\begin{tabular}{|C{3cm}|C{2.5cm}C{2.5cm}||C{4cm}|}
  \hline
   & \multicolumn{2}{c||}{Parameter-estimations} & \\
  \hline
 Variables & ${d^1}'$ & ${a^1}'$ & Correlations with $\widetilde{f^1}$\\
  \hline
 altitude & 4.43 & 0.62 & 0.66\\
 pluvio\_yr & 44.45 & 0.16 & 0.17 \\
 pluvio\_1 & 2.48 & -0.91 & -0.97 \\
 pluvio\_2 & 4.32 & -0.88 & -0.94 \\
 pluvio\_3 & 9.65 & -0.47 & -0.5 \\
 pluvio\_4 & 8.56 &  -0.28 & -0.3 \\
 pluvio\_5 & 6.68 & 0.26 & 0.28 \\
 pluvio\_6 & 5.98 & 0.83 & 0.89 \\
 pluvio\_7 & 4.78  & 0.81 & 0.86 \\
 pluvio\_8 & 4.17 & 0.86 & 0.91 \\
 pluvio\_9 & 11.46 & 0.72 & 0.77 \\
 pluvio\_10 & 10.17 & 0.34 & 0.36 \\
 pluvio\_11 & 4.36 & -0.83 & -0.88 \\
 pluvio\_12 & 2.13 & -0.9 & -0.96 \\
 lon & 14.57 & 0.04 & 0.04 \\
 lat & 2.49 & 0.92 & 0.98 \\
  \hline  
\end{tabular}
\captionof{table}{Application to the \textit{genus} data with geologic covariate : estimations of parameters ${d^1}'$ and ${a^1}'$, and correlations between $\widetilde{f^1}$ and the variables $X^1$}
\label{tabl_estim_param_d1a1_avec}

\begin{figure}[H]
\centering
\includegraphics[width = .8\textwidth]{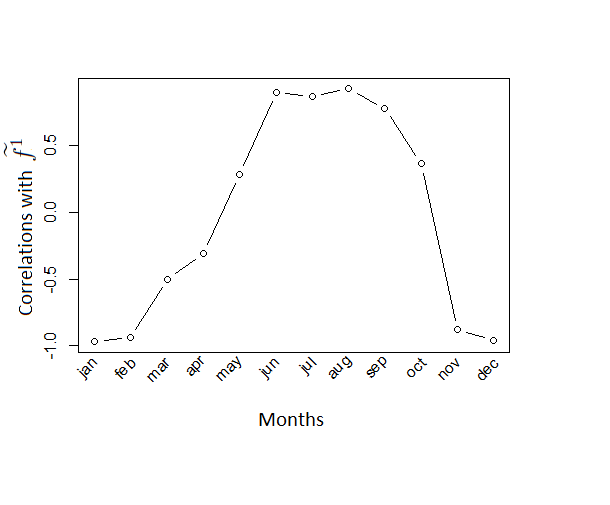}  
\caption{Correlations between $\widetilde{f^1}$ and the monthly variables of $X^1$ : two rainfall regimes.}
\label{Corr_f1_tilde}
\end{figure}

\begin{tabular}{|C{1cm}C{1cm}
C{.8cm}C{.8cm}C{.8cm}|}
\hline
 \multicolumn{5}{|c|}{Scalar parameter-estimations}\\
   \hline
 $c^1$ & $c^2$ & $\sigma^2_{1}$ &  $\sigma^2_{2}$ & $\sigma^2_{Y}$  \\
  \hline
 0.35 & 0.01 & 0.50 & 0.53 & 0.84 \\
  \hline
\end{tabular}
\captionof{table}{Application to \textit{genus} data with geologic covariate: scalar parameter-estimations }
\label{tabl_estim_param_scal_avec}
\end{center}

It can be seen in Tables \ref{tabl_estim_param_db_avec} and \ref{tabl_etendue_avec} that for certain species, the geologic substrate seems to be of great importance (e.g. for gen1, gen5, gen7, gen9, gen12, gen16, gen21, gen25, gen26, gen27), whereas for others, it only has a small impact on the abundances (e.g. for gen2, gen6, gen8, gen10, gen18, gen20, gen23). 
Moreover, Table \ref{tabl_estim_param_db_avec} shows that the correlations between $\widetilde{g}$ and $Y$ are high in absolute value only for few variables : \textit{gen2}, \textit{gen23}, \textit{gen24} and \textit{gen25}. Therefore, only these are well accounted for by our model. Although we have carried out the analysis with variables \textit{gen2}, \textit{gen3}, \textit{gen8}, \textit{gen10}, \textit{gen11}, \textit{gen15}, \textit{gen17}, \textit{gen23}, \textit{gen24} and \textit{gen25}, the results are practically the same when we take all variables.
In table \ref{tabl_estim_param_d1a1_avec}, the correlations between $\widetilde{f^1}$ and variables \textit{pluvio\_1} to \textit{pluvio\_12} of $X^1$ show two rainfall regimes. Indeed, \textit{pluvio\_1} corresponds to january, \textit{pluvio\_2}, to february, ... \textit{pluvio\_12} to december. The Central African Republic has a tropical climate : the dry season ranges from November to April and the rainy season from June to September. Figure \ref{Corr_f1_tilde} shows  that $\widetilde{f^1}$ is positively correlated to the rainfalls of the rainy season and negatively to those of the dry one.

\subsection{Model without covariate}
What if we omit the geologic substrate as covariate? 
\subsubsection{Model specification}
We now consider the model without covariate \textit{geology}, i.e. every $T$ is reduced to $\mathds {1}_n$:
\begin{equation}
\left\{ \begin{array}{ccl}
Y &=& \mathds {1}_n d' + g b' + \varepsilon^{Y} \nonumber\\
X^1 &=& \mathds {1}_n {d^1}' + f^1 {a^1}' + \varepsilon^{1} \nonumber\\ 
X^2 &=& \mathds {1}_n {d^2}' + f^{2} {a^2}' + \varepsilon^{2} \nonumber\\
g &=& f^1 c^1 + f^2 c^2 + \varepsilon^{g} \nonumber
\end{array}
\right.
\end{equation}

Where, there always are $n=1000$, $q_Y=27$, $q_1=16$, $q_2=23$. Also, $d'$ (resp. ${d^1}'$ and ${d^2}'$) are parameters vectors that contains the means of the $Y$'s (resp. $X^1$ and $X^2$). 

\subsubsection{Results}
With a threshold value $\varepsilon=10^{-3}$, convergence was reached after 49 iterations. 
Some parameter-estimations are presented in Table \ref{tabl_estim_param_db_sans} and \ref{tabl_estim_param_scal_sans}. For practical reasons, in Appendix D tables \ref{tabl_estim_param_d1a1_sans} and \ref{tabl_estim_param_d2a2_sans} present the remainder of the parameter-estimations.

\begin{center}
\begin{tabular}{|C{3cm}|C{2.5cm}C{2.5cm}||C{4cm}|}
  \hline
   & \multicolumn{2}{c||}{Parameter-estimations} & \\
   \hline
 Variables & $d'$ & $b'$ & Correlation with $\widetilde{g}$\\
  \hline
  gen1 & 0.95 & 0.16 & 0.19\\
  gen2 & 0.45 & 0.34 & 0.41\\
  gen3 &  0.38 & 0.34 & 0.42\\
  gen4 &  0.35 & 0.18 & 0.22\\
  gen5 &  0.45 & -0.15 & -0.18\\
  gen6 &  0.54 & 0.18 & 0.22\\
  gen7 & 0.44 & -0.16 & -0.2\\
  gen8 &  0.51 & 0.29 & 0.35\\
  gen9 & 0.56 & -0.01 & -0.01\\
  gen10 & 0.86 & -0.2 & -0.24\\
  gen11 & 1.54 & 0.35 & 0.43\\
  gen12 & 0.74 & -0.04 & -0.04\\
  gen13 & 0.94 & 0.28 & 0.35\\
  gen14 & 0.35 & 0.1 & 0.12\\
  gen15 & 0.54 & 0.39 & 0.47\\
  gen16 &  0.58 & 0.32 & 0.4\\
  gen17 & 0.48 & 0.17 & 0.21\\
  gen18 & 0.17 & 0.16 & 0.2\\
  gen19 & 0.41 & 0.04 & 0.05\\
  gen20 & 0.53 & -0.13 & -0.16\\
  gen21 & 0.93 & 0.41 & 0.5\\
  gen22 & 0.40 & 0.39 & 0.48\\
  gen23 & 0.96 & 0.45 & 0.56\\
  gen24 & 0.74 & 0.54 & 0.66\\
  gen25 & 0.76 & 0.64 & 0.78\\
  gen26 & 0.85 & -0.33 & -0.41\\
  gen27 & 0.40 & -0.19 & -0.24\\
 \hline
\end{tabular}
\captionof{table}{Application to the \textit{genus} data without covariate : estimations of parameters $D'$ and $b'$, and correlations between $\widetilde{g}$ and the variables $Y$}
\label{tabl_estim_param_db_sans} 
\end{center}
\begin{center}
\begin{tabular}{|C{1cm}C{1cm}
C{.8cm}C{.8cm}C{.8cm}|}
\hline
 \multicolumn{5}{|c|}{Scalar parameter-estimations}\\
   \hline
 $c^1$ & $c^2$ & $\sigma^2_{1}$ &  $\sigma^2_{2}$ & $\sigma^2_{Y}$  \\
  \hline
 0.27 & -0.07 & 0.50 & 0.53 & 0.90 \\
  \hline
\end{tabular}
\captionof{table}{Application to the \textit{genus} data without covariate : scalar estimations of parameters}
\label{tabl_estim_param_scal_sans}
\end{center}

Table \ref{tabl_estim_param_scal_avec} and \ref{tabl_estim_param_scal_sans} show that, the geological effect is being considered or removed, the geographic factor (position and rainfalls) keeps a much greater effect than the EVI's.
The estimations of $\sigma^2_{1}$,  $\sigma^2_{2}$, $\sigma^2_{Y}$ don't change significantly. For both, the rainfall regimes are identically identified. However, table \ref{tabl_estim_param_d1a1_avec} shows the great impact of geologic substrate on the abundance of species gen1, gen5, gen7, gen9, gen12, gen16, gen21, gen25, gen26 and gen27. Therefore, the presence of covariate \textit{geology} in the model is relevant. 

\subsection{Assessing the model quality through re-sampling}
To assess the stability of results and thus, validate the models (with covariate), we use a 5-fold re-sampling technique: 5 separate samples are randomly extracted from the complete \textit{genus} data, thus, their size is kept equal to 200 units. For each of them, we obtain estimated parameters and factors. Then, for each sample, we compute an average Mean Square Error (MSE) and an average correlation between the parameter-estimates obtained on the sample and those obtained on the complete data. Finally, on each sample, we calculate an average MSE and correlation between the factor-estimates obtained on the sample and the corresponding ones obtained on the complete data for the units belonging to the sample.

\subsubsection{Model with geologic covariate}
\begin{figure}[H]
\centering 
\includegraphics[width = .5\textwidth]{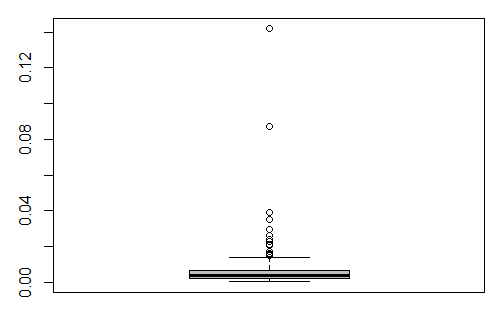}  
\caption{Box plot of the average MSE's between the parameter-estimates obtained on the 5 \textit{genus} data sub-samples and those obtained on the complete data.}
\label{MSE_boxplot_theta_ac_cov}
\end{figure}

\begin{figure}[H]
\centering
\includegraphics[width = .5\textwidth]{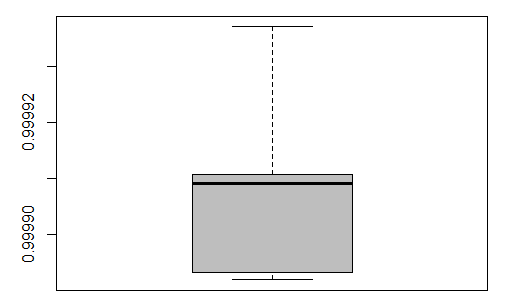}  
\caption{Box plot of the average correlations between the parameter-estimates obtained on the 5 \textit{genus} data sub-samples and those obtained on the complete data.}
\label{Corr_boxplot_theta_ac_cov}
\end{figure}

Figure \ref{MSE_boxplot_theta_ac_cov} (resp. \ref{Corr_boxplot_theta_ac_cov}) shows the average 
MSE (resp. the correlation) between estimated parameters on 5 data samples $\theta^*_{s \in \llbracket 1,5 \rrbracket}$ and estimated parameters on the complete data $\theta^*$. More precisely, for these average MSE (respectively correlation), the median is $3.85*10^{-3}$ (resp. $0.99$), the first quartile is $1.95*10^{-3}$ (resp. $0.99$) and the third quartile is $6.17*10^{-3}$ (resp. $0.99$). These values are close to $0$ (resp. $1$).  So, we can be rather confident in the estimates of parameters obtained in the last section. 

\begin{figure}[H]
\centering
\includegraphics[width = .5\textwidth]{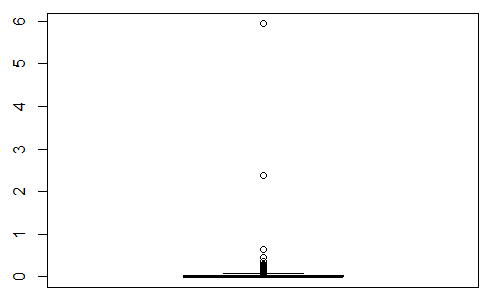}
\caption{Box-plot of the average MSE of factor-estimates.} 
\label{MSE_boxplot_factors_ac_cov}
\end{figure}

\begin{figure}[H]
\centering 
\includegraphics[width = .5\textwidth]{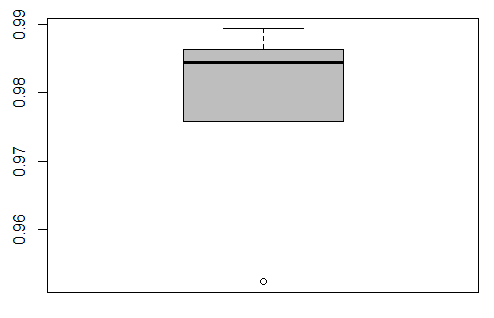}
\caption{Box-plot of the average correlation of factor-estimates.}
\label{Corr_boxplot_factors_ac_cov}
\end{figure}

 Figure \ref{MSE_boxplot_factors_ac_cov} and Figure \ref{Corr_boxplot_factors_ac_cov} respectively give the box-plot of the factors' average MSE and correlation for each of the 5 samples. More precisely, for these average MSE's (respectively correlations), the median is $1.15*10^{-2}$ (resp. $0.98$), the first quartile is $7.44*10^{-3}$ (resp. $0.98$) and the third quartile is $3.53*10^{-2}$ (resp. $0.99$). These values are close enough to $0$ (resp. $1$) to allow us to be confident in the estimates obtained on the complete data.


\section{Conclusion}
The maximum-likelihood estimation method is known to be a stringent method of estimation having nice properties. In the context of estimation methods of SEM, the LISREL approach is based on likelihood maximization, contrary to PLS, THEME, and other component-based methods. However, LISREL only focusses on the variance-covariance structure of the data and does not allow to estimate the LV's, contrary to PLS and THEME. To estimate them, we proposed to carry out likelihood maximization through the EM algorithm. This approach assumes that VLs are factors, which is constraining than assuming they are components. Therefore this approach combines the stringency of likelihood maximization with the possibility to estimate the LV's. This presented this new approach and performed sensitivity analysis to show its performances. Then, an application on environmental data was made, which shows how to use this method. Along with re-sampling for validation purposes. 

\section{Acknowledgments}
The forest inventories were funded by the French Agency for Development (AFD) through the PARPAF project. We wish to
thank S. Chong (TCA), A. Banos (SCAF), and the ‘"Ministère des Eaux, Forêts, Chasse et Pêche" of the Central African Republic
for authorizing access to the inventory data, and the field teams who drew up these inventories. This study is part of the ErA
Net BiodivERsA CoForChange project, funded by the National Research Agency (ANR) and the Natural Environment Research
Council (NERC), involving 16 European, African and international partners and a number of timber companies (see the list
on the website, http://www.coforchange.eu).


\section*{Appendix A. Calculation of the complete data log-likelihood function ${\cal L}$}
\begin{proof}
In our case $p=2$, $\psi_Y=\sigma^2_Y Id_{q_Y}$, $\psi_{1}=\sigma^2_1 Id_{q_1}$ and $\psi_{2}=\sigma^2_2 Id_{q_2}$, and for observation $i$, the model is formulated as follows:
\begin{equation}
\left\{ \begin{array}{ccl}
y_i' &=& {t_i}'D+ g_i b' + {{\varepsilon_i}^{y}}'  \nonumber\\
{x_i^1}' &=& {t_i^1}'D^1 + f_i^1 {a^1}' + {{\varepsilon_i}^{1}}'  \nonumber\\
{x_i^2}' &=& {t_i^2}'D^2  + f_i^2 {a^2}' + {{\varepsilon_i}^{2}}'  \nonumber\\
g_i &=& f_i^1 c^1 + f_i^2 c^2 + {{\varepsilon_i}^{g}}  \nonumber
\end{array}
\right.
\end{equation}

We have,
\begin{align*}
\begin{split}
p(z_i,h_i;\theta) 
&= p(y_i,x_i^1,x_i^2,g_i,f_i^1,f_i^2;\theta)\\
&= p(y_i,x_i^1,x_i^2|g_i,f_i^1,f_i^2;\theta)p(g_i,f_i^1,f_i^2;\theta)\\
&= p(y_i,x_i^1,x_i^2|g_i,f_i^1,f_i^2;\theta)p(g_i|f_i^1,f_i^2;\theta)p(f_i^1,f_i^2;\theta)\\
&= p(y_i,x_i^1,x_i^2|g_i,f_i^1,f_i^2;\theta)p(g_i|f_i^1,f_i^2;\theta)p(f_i^1;\theta)p(f_i^2;\theta)\\
&= p(y_i,x_i^1,x_i^2|g_i,f_i^1,f_i^2;\theta)p(g_i|f_i^1,f_i^2;\theta)p(f_i^1)p(f_i^2)\\
&= p(x_i^1,x_i^2|y_i,g_i,f_i^1,f_i^2;\theta)p(y_i|g_i,f_i^1,f_i^2;\theta)p(g_i|f_i^1,f_i^2;\theta)p(f_i^1)p(f_i^2)\\
&= p(x_i^1,x_i^2|f_i^1,f_i^2;\theta)p(y_i|g_i;\theta)p(g_i|f_i^1,f_i^2;\theta)p(f_i^1)p(f_i^2)\\
&= p(x_i^1|x_i^2,f_i^1,f_i^2;\theta)p(x_i^2|f_i^1,f_i^2;\theta)p(y_i|g_i;\theta)p(g_i|f_i^1,f_i^2;\theta)p(f_i^1)p(f_i^2)\\
&= p(x_i^1|f_i^1;\theta)p(x_i^2|f_i^2;\theta)p(y_i|g_i;\theta)p(g_i|f_i^1,f_i^2;\theta)p(f_i^1)p(f_i^2)\\
\end{split}
\end{align*}
Where $\theta=\lbrace D, D^1, D^2,  b, {a^1}, {a^2}, {c^1}, {c^2}, \psi_Y, \psi_{1}, \psi_{2} \rbrace$ is the set of model parameters.
Therefore,
\begin{align*}
\begin{split}
&{\cal L}(\theta;z_i,h_i) = {\cal L}(\theta;x_i^1|f_i^1)+{\cal L}(\theta;x_i^2|f_i^2)+{\cal L}(\theta;y_i|g_i)+{\cal L}(\theta;g_i|f_i^1,f_i^2)+{\cal L}(f_i^1)+{\cal L}(f_i^2)\\
\end{split}
\end{align*}
Because of the model and the normal distribution properties we obtain:\\
$x_i^m|f_i^m \sim {\cal N}({t_i^m}'D^m + f_i^m {a^m}', \psi_{X^m})$\\ 
$y_i|g_i \sim {\cal N}({t_i}'D+ g_i b', \psi_Y)$\\ 
$g_i|f_i^1,f_i^2 \sim {\cal N}( f_i^1 c^1 + f_i^2 c^2 , 1)$\\
$f_i^m \sim {\cal N}(0, 1)$\\
Then, we obtain the complete data log-likelihood function: 
\begin{align*}
\begin{split}
{\cal L}(\theta;z,h) = - \frac{1}{2} 
\sum \limits_{\substack{i=1}}^{n}
&\lbrace
ln|\psi_Y| + ln|\psi_{1}| + ln|\psi_{2}| \\
&+ (y_i - D't_i - g_ib)'\psi_Y^{-1}(y_i - D't_i - g_ib) \\
&+ (x_i^1 - {D^1}'t_i^1 - f_i^1a^1)'\psi_{1}^{-1}(x_i^1 - {D^1}'t_i^1 - f_i^1a^1) \\
&+ (x_i^2 - {D^2}'t_i^2 - f_i^2a^2)'\psi_{2}^{-1}(x_i^2 - {D^2}'t_i^2 - f_i^2a^2) \\
&+ (g_i - c^1\ f_i^1 - c^2\ f_i^2)^2 + (f_i^1)^2 + (f_i^2)^2 \rbrace + \lambda
\\
\end{split}
\end{align*}
Where $\lambda$ a constant. Also, the set of model parameters\\ $\theta=\lbrace D, D^1, D^2,  b, {a^1}, {a^2}, {c^1}, {c^2}, \psi_Y, \psi_{1}, \psi_{2} \rbrace$ in our case corresponds to\\ $\theta=\lbrace D, D^1, D^2,  b, {a^1}, {a^2}, {c^1}, {c^2},\sigma^2_Y, \sigma^2_1, \sigma^2_2 \rbrace$ because of the simplification made in the section 2.3, in our case. Indeed, $\psi_Y=\sigma^2_Y Id_{q_Y}$, $\psi_{1}=\sigma^2_1 Id_{q_1}$ and $\psi_{2}=\sigma^2_2 Id_{q_2}$.\\
Therefore, we can also write the  complete data log-likelihood function:\\
\begin{align*}
\begin{split}
{\cal L}(\theta;z,h) = - \frac{1}{2} 
\sum \limits_{\substack{i=1}}^{n}
&\lbrace
q_Yln(\sigma^2_Y) + q_1ln(\sigma^2_1) + q_2ln(\sigma^2_2) \\
&+ {\sigma^{-2}}_Y(y_i - D't_i - g_ib)'(y_i - D't_i - g_ib) \\
&+ {\sigma^{-2}}_1(x_i^1 - {D^1}'t_i^1 - f_i^1a^1)'(x_i^1 - {D^1}'t_i^1 - f_i^1a^1) \\
&+ {\sigma^{-2}}_2(x_i^2 - {D^2}'t_i^2 - f_i^2a^2)'(x_i^2 - {D^2}'t_i^2 - f_i^2a^2) \\
&+ (g_i - c^1\ f_i^1 - c^2\ f_i^2)^2 + (f_i^1)^2 + (f_i^2)^2 \rbrace + \lambda
\\
\end{split}
\end{align*}
\end{proof}

\section*{Appendix B. Demonstration of the normality of the distribution of $h_i|z_i$}
\begin{proof}
In our case $p=2$, $\psi_Y=\sigma^2_Y Id_{q_Y}$, $\psi_{1}=\sigma^2_1 Id_{q_1}$ and $\psi_{2}=\sigma^2_2 Id_{q_2}$, and for observation $i$, the model is formulated as follows:
\begin{equation}
\left\{ \begin{array}{ccl}
y_i' &=& {t_i}'D+ g_i b' + {{\varepsilon_i}^{y}}'  \nonumber\\
{x_i^1}' &=& {t_i^1}'D^1 + f_i^1 {a^1}' + {{\varepsilon_i}^{1}}'  \nonumber\\
{x_i^2}' &=& {t_i^2}'D^2  + f_i^2 {a^2}' + {{\varepsilon_i}^{2}}'  \nonumber\\
g_i &=& f_i^1 c^1 + f_i^2 c^2 + {{\varepsilon_i}^{g}}  \nonumber
\end{array}
\right.
\end{equation}

To prove the normality of the distribution of $h_i|z_i$ presented in section 3.1.2., we use the classical result\footnote{If two variables $X_1$ and $X_2$ are normally distributed such that, 

$\begin{pmatrix} X_1 \\ X_2  \end{pmatrix} \sim {\cal N}(\mu = \begin{pmatrix} \mu_1 \\ \mu_2  \end{pmatrix},\Sigma = \begin{pmatrix} \Sigma_{11} & \Sigma_{12} \\ \Sigma_{21} & \Sigma_{22} \end{pmatrix})$\\
where, $\mu_1$ ($r \times 1$), $\mu_2$ ($s \times 1$), $\Sigma_{11}$ ($r \times r$), $\Sigma_{12}$($r \times s$), $\Sigma_{21}$ ($s \times r$) and $\Sigma_{22}$ ($s \times s$);\\
then, 
\begin{equation}
\hspace{-.15cm}
(X_1|X_2=x_2)\sim {\cal N}(M=\mu_1+\Sigma_{12} {\Sigma_{22}}^{-1}(x_2-\mu_2), \phi=\Sigma_{11}-\Sigma_{12} {\Sigma_{22}}^{-1}\Sigma_{21})
\label{dist_cond}
\end{equation}
} about the conditioning of normally distributed variables.
Before using this result, we calculate the joint distribution of $(g_i, f^1_i, f^2_i,y_i,x^1_i, x^2_i)$.\\ 

We know that, for observation $i$,\\
$y_i \sim {\cal N}(D't_i , b((c^1)^2 + (c^2)^2 + 1)b' + \Psi_{Y})$\\
$x^m_i \sim {\cal N}({D^m}'t^m_i , a^m{a^m}' + \Psi_{m})$\\
$g_i \sim {\cal N}( 0 ,  (c^1)^2 + (c^2)^2 + 1)$\\
$f^m_i \sim {\cal N}( 0 , 1)$\\
Then, after compute the required covariances we obtain, \\
$(g_i, f^1_i, f^2_i) \sim {\cal N}( \begin{pmatrix} 0 \\ 0 \\ 0  \end{pmatrix}, \begin{pmatrix}(c^1)^2 + (c^2)^2 + 1 & c^1 & c^2 \\ c^1 & 1 & 0 \\  c^2 & 0 & 1 \end{pmatrix})$\\
And,\\
$(y_i,x^1_i, x^2_i) \sim {\cal N}( \begin{pmatrix} D't_i \\ {D^1}'t^1_i \\ {D^2}'t^2_i\end{pmatrix}, \begin{pmatrix} ((c^1)^2 + (c^2)^2 + 1)bb' + \Psi_{Y} & c^1b{a^1}' & c^2b{a^2}' \\ c^1a^1b' & a^1{a^1}' + \Psi_{1} & 0_{(q_1,q_2)} \\  c^2a^2b' & 0_{(q_2,q_1)} & a^2{a^2}' + \Psi_{2} \end{pmatrix})$\\

Then, after compute the required covariances we obtain the joint distribution,
$(g_i, f^1_i, f^2_i, y_i,x^1_i, x^2_i) \sim {\cal N}(M^*_i, \Sigma^*
)$ such as,\\
$M^*_i = \begin{pmatrix} 0_{(3,1)} \\ D't_i \\ {D^1}'t^1_i \\ {D^2}'t^2_i \end{pmatrix}$
and
$\Sigma^* = \begin{pmatrix}
\Sigma^*_1 & \Sigma^*_2 \\ 
{\Sigma^*_2}' & \Sigma^*_3
  \end{pmatrix}$. Where,\\
$\Sigma^*_1 = \begin{pmatrix} (c^1)^2+ (c^2)^2+1 & c^1 & c^2 \\  
                                 c^1 & 1 & 0 \\
                                 c^2 & 0 & 1
                                 \end{pmatrix}$\\
$\Sigma^*_2 = \begin{pmatrix} ((c^1)^2+ (c^2)^2 + 1)b'  &  c^1{a^1}' & c^2{a^2}' \\  
                        c^1b' &  {a^1}' & 0_{(1,q_2)}\\
                        c^2b' & 0_{(1,q_1)} & {a^2}' \end{pmatrix}$\\
$\Sigma^*_3 = \begin{pmatrix} ((c^1)^2+ (c^2)^2 + 1)bb' + \Psi_{Y} & c^1b{a^1}' & c^2b{a^2}' \\
                       c^1a^1b' & a^1{a^1}' + \Psi_{1} &  0_{(q_1,q_2)}\\
                      c^2a^2b' & 0_{(q_2,q_1)} & a^2{a^2}' + \Psi_{2}
                       \end{pmatrix}$\\                       
                                 
Finally, we use result \eqref{dist_cond} and obtain the distribution, $h_i|z_i\sim {\cal N}(M_i,\Sigma)$ where, $ M_i = \Sigma^*_2 {\Sigma^*_3}^{-1} \mu^*_i$ and $\Sigma = \Sigma^*_1 - \Sigma^*_2 {\Sigma^*_3}^{-1}{\Sigma^*_2}'$, such that
$\mu^*_i = \begin{pmatrix} y_i - D't_i  \\
                       x_i^1 - {D^1}'t_i^1 \\
                       x_i^2 - {D^2}'t_i^2 \end{pmatrix}$\\

\end{proof}

\section*{Appendix C. Calculation of the first-order derivatives of ${\cal L}$}
\begin{proof}
We search the first-order derivatives of the complete data log-likelihood function: 
\begin{align*}
\begin{split}
{\cal L}(\theta;z,h) = - \frac{1}{2} 
\sum \limits_{\substack{i=1}}^{n}
&\lbrace
ln|\psi_Y| + ln|\psi_{1}| + ln|\psi_{2}| \\
&+ (y_i - D't_i - g_ib)'\psi_Y^{-1}(y_i - D't_i - g_ib) \\
&+ (x_i^1 - {D^1}'t_i^1 - f_i^1a^1)'\psi_{1}^{-1}(x_i^1 - {D^1}'t_i^1 - f_i^1a^1) \\
&+ (x_i^2 - {D^2}'t_i^2 - f_i^2a^2)'\psi_{2}^{-1}(x_i^2 - {D^2}'t_i^2 - f_i^2a^2) \\
&+ (g_i - c^1\ f_i^1 - c^2\ f_i^2)^2 + (f_i^1)^2 + (f_i^2)^2 \rbrace + \lambda
\\
\end{split}
\end{align*}
Where $\lambda$ constant, $\theta=\lbrace D, D^1, D^2,  b, {a^1}, {a^2}, {c^1}, {c^2}, \psi_Y, \psi_{1}, \psi_{2} \rbrace$, $\psi_Y=\sigma^2_Y Id_{q_Y}$, $\psi_{1}=\sigma^2_1 Id_{q_1}$ and $\psi_{2}=\sigma^2_2 Id_{q_2}$.\\
Therefore, there are matrix-parameters ($D, D^1, D^2$), vector-parameters ($b, {a^1}, {a^2}$) and scalar parameters (${c^1}, {c^2}, \sigma^2_Y, \sigma^2_1, \sigma^2_2$). Then, ${\cal L}$ is a sum of three types of functions: the logarithm, the square function and a quadratic form function $(w - X\beta)'\Gamma(w - X\beta)$, where $\Gamma$ is symmetric and $w$ ($q \times 1$), $X$ ($q \times m$), $\beta$ ($m \times 1$) and $\Gamma$ ($q \times q$). The first-order derivatives of the logarithm function  and the square function are in our case trivial. The first-order derivative of $(w - X\beta)'\Gamma(w - X\beta)$ by $X$ is less trivial but necessary. Let us start by making explicit the first-order derivative of $(w - X\beta)'\Gamma(w - X\beta)$ with respect to $X$.
\begin{align*}
\begin{split}
d_X[(w - X\beta)'\Gamma(w - X\beta)] 
&= (w - X\beta)'\Gamma(- dX\beta) + (- dX\beta)'\Gamma(w - X\beta)\\
&= -2(w - X\beta)'\Gamma(dX\beta) \\
&= tr[-2(w - X\beta)'\Gamma(dX\beta)] \\
&= tr[-2\beta(w - X\beta)'\Gamma dX] \\
&= <-2\beta(w - X\beta)'\Gamma | dX>
\end{split}
\end{align*}
Therefore,
\begin{align*}
\begin{split}
\frac{d}{dX}[(w - X\beta)'\Gamma(w - X\beta)] 
&= (-2\beta(w - X\beta)'\Gamma)'\\
&= -2(\beta(w - X\beta)'\Gamma)' \\
&= -2\Gamma(w - X\beta){\beta}'
\end{split}
\end{align*}
Likewise, we establish that :
\begin{equation*}
\frac{\partial}{\partial D'}{\cal L}(z,h) = \sum \limits_{\substack{i=1}}^{n}
\psi_Y^{-1}(y_i - D't_i - g_ib) {t_i}'
\end{equation*}

Similar reasoning can be applied to $D^m$ and allows to obtain the second row of \eqref{syst_derrivees1}. Concerning the third and the fourth row of \eqref{syst_derrivees1}, we use the classical result : 
\begin{equation*}
\frac{\partial}{\partial \beta}[(w-X\beta)'\Gamma(w-X\beta)] = -2X'\Gamma(w-X\beta)
\end{equation*}
Eventually, the fifth, the sixth and the eighth rows of \eqref{syst_derrivees1} are obtained in a trivial way.
\end{proof}

\section*{Appendix D. Tables of section 5.}

\begin{center}
\begin{tabular}{|C{1.5cm}|C{3cm}|}
   \hline
 Variables & Differences  \\
  \hline
  gen1 & \colorbox{gris25}{0.80}   \\
  gen2 & \textit{0.28}\\
  gen3 & 0.62  \\
  gen4 & 0.52 \\
  gen5 & \colorbox{gris25}{1.04} \\
  gen6 & \textit{0.32} \\
  gen7 & \colorbox{gris25}{0.80}\\
  gen8 & \textit{0.20} \\
  gen9 & \colorbox{gris25}{0.92} \\
  gen10 & \textit{0.40} \\
  gen11 & 0.67\\
  gen12 & \colorbox{gris25}{1.24} \\
  gen13 & 0.53 \\
  gen14 & 0.48 \\
  gen15 & 0.45 \\
  gen16 & \colorbox{gris25}{0.82} \\
  gen17 & 0.54\\
  gen18 & \textit{0.25} \\
  gen19 & 0.52 \\
  gen20 & \textit{0.35} \\
  gen21 & \colorbox{gris25}{1.01}\\
  gen22 & 0.49\\
  gen23 & \textit{0.31} \\
  gen24 & 0.69 \\
  gen25 & \colorbox{gris25}{1.46}  \\
  gen26 & \colorbox{gris25}{1.52}  \\
  gen27 & \colorbox{gris25}{0.93} \\
 \hline
\end{tabular}
\captionof{table}{Application to the \textit{genus} data with geologic covariate : Differences between maximal and minimal values of geologic effects $D[1,]$, $D[1,]+D[2,]$, $D[1,]+D[3,]$, $D[1,]+D[4,]$, $D[1,]+D[5,]$ (highlights on the greater differences, italics on the smaller)}
\label{tabl_etendue_avec}
\end{center}

\begin{tabular}{|C{3cm}|C{2.5cm}C{2.5cm}||C{4cm}|}
  \hline
   & \multicolumn{2}{c||}{Parameter-estimations} &  \\
  \hline
 Variables & ${d^2}'$ & ${a^2}'$ & Correlations with $\widetilde{f^2}$\\
  \hline
evi\_1 & 15.51 & 0.63 & 0.65\\
evi\_2 & 13.47 & 0.59 & 0.6\\
evi\_3 & 14.83 & 0.51 & 0.52\\
evi\_4 & 14.67 & 0.58 & 0.6\\
evi\_5 & 16.44 & 0.56 & 0.57\\
evi\_6 & 18.74 & 0.51 & 0.52\\
evi\_7 & 18.44 & 0.75 & 0.76\\
evi\_8 & 20.59 & 0.8 &  0.82\\
evi\_9 & 21.83 & 0.76 & 0.78\\
evi\_10 & 19.19 & 0.74 & 0.76\\
evi\_11 & 18.22 & 0.67 & 0.69\\
evi\_12 & 15.92 & 0.61 & 0.63\\
evi\_13 & 15.4 & 0.58 & 0.6\\
evi\_14 & 13.51 &  0.7 & 0.72\\
evi\_15 & 14.57 & 0.69 & 0.71\\
evi\_16 & 14.95 & 0.76 & 0.78\\
evi\_17 & 16.09 & 0.73 & 0.75\\
evi\_18 & 15.95 & 0.77 & 0.79\\
evi\_19 & 17.12 & 0.73 & 0.75\\
evi\_20 & 15.02 & 0.75 & 0.77\\
evi\_21 & 15.87 & 0.75 & 0.77\\
evi\_22 & 14.21 & 0.71 & 0.73\\
evi\_23 & 15.26 & 0.68 & 0.69\\
  \hline
\end{tabular}
\captionof{table}{Application to the \textit{genus} data with geologic covariate : estimations of parameters ${d^2}'$ and ${a^2}'$, and correlations between $\widetilde{f^2}$ and the variables $X^2$}
\label{tabl_estim_param_d2a2_avec}

\begin{tabular}{|C{3cm}|C{2.5cm}C{2.5cm}||C{4cm}|}
  \hline
   & \multicolumn{2}{c|}{Parameter-estimations} & \\
  \hline
 Variables & ${d^1}'$ & ${a^1}'$ & Correlations with $\widetilde{f^1}$\\
  \hline
 altitude & 4.43 & 0.63 & 0.66 \\
 pluvio\_yr & 44.45 & 0.16 & 0.17\\
 pluvio\_1 & 2.48 & -0.92 & -0.97\\
 pluvio\_2 & 4.32 & -0.89 & -0.94\\
 pluvio\_3 & 9.65 & -0.48 & -0.5\\
 pluvio\_4 & 8.56 &  -0.29 & -0.31\\
 pluvio\_5 & 6.68 & 0.26 & 0.28\\
 pluvio\_6 & 5.98 & 0.84 & 0.89\\
 pluvio\_7 & 4.78  & 0.82 & 0.86\\
 pluvio\_8 & 4.17 & 0.87 & 0.92\\
 pluvio\_9 & 11.46 & 0.73 & 0.77\\
 pluvio\_10 & 10.17 & 0.34 & 0.36\\
 pluvio\_11 & 4.36 & -0.84 & -0.88\\
 pluvio\_12 & 2.13 & -0.91 & -0.96\\
 lon & 14.57 & 0.04 & 0.05\\
 lat & 2.49 & 0.93 & 0.98\\
\hline  
\end{tabular}
\captionof{table}{Application to the \textit{genus} data without covariate : estimations of parameters ${d^1}'$ and ${a^1}'$, and correlations between $\widetilde{f^1}$ and the variables $X^1$}
\label{tabl_estim_param_d1a1_sans}

\begin{tabular}{|C{3cm}|C{2.5cm}C{2.5cm}||C{4cm}|}
  \hline
   & \multicolumn{2}{c||}{Parameter-estimations} & \\
  \hline
 Variables & ${d^2}'$ & ${a^2}'$ & Correlations with $\widetilde{f^2}$\\
  \hline
evi\_1 & 15.51 & 0.63 & 0.65\\
evi\_2 & 13.47 & 0.59 & 0.6\\
evi\_3 & 14.83 & 0.51 & 0.52\\
evi\_4 & 14.67 & 0.58 & 0.6\\
evi\_5 & 16.44 & 0.56 & 0.57\\
evi\_6 & 18.74 & 0.5 & 0.52\\
evi\_7 & 18.44 & 0.74 & 0.76\\
evi\_8 & 20.59 & 0.8 & 0.82\\
evi\_9 & 21.83 & 0.76 & 0.78\\
evi\_10 & 19.19 & 0.74 & 0.76\\
evi\_11 & 18.22 & 0.67 & 0.69\\
evi\_12 & 15.92 & 0.61 & 0.63\\
evi\_13 & 15.4 & 0.58 & 0.6\\
evi\_14 & 13.51 &  0.7 & 0.72\\
evi\_15 & 14.57 & 0.69 & 0.71\\
evi\_16 & 14.95 & 0.76 & 0.78\\
evi\_17 & 16.09 & 0.73 & 0.75\\
evi\_18 & 15.95 & 0.77 & 0.79\\
evi\_19 & 17.12 & 0.73 & 0.75\\
evi\_20 & 15.02 & 0.75 & 0.77\\
evi\_21 & 15.87 & 0.75 & 0.77\\
evi\_22 & 14.21 & 0.71 & 0.73\\
evi\_23 & 15.26 & 0.67 & 0.69\\
  \hline
\end{tabular}
\captionof{table}{Application to the \textit{genus} data without geologic covariate : estimations of parameters ${d^2}'$ and ${a^2}'$, and correlations between $\widetilde{f^2}$ and the variables $X^2$}
\label{tabl_estim_param_d2a2_sans}


\section{References}
\nocite{*}
{\footnotesize
\bibliography{Article_SEM_VCX_complet}
}
\bibliographystyle{elsarticle-harv} 

%
%
%


\end{document}